\documentclass[12pt]{article}
\usepackage{amsthm, amsmath, amssymb}
\newtheorem{theorem}{Theorem}[section]
\newtheorem{lemma}[theorem]{Lemma}
\newtheorem{corollary}[theorem]{Corollary}
\newtheorem{proposition}[theorem]{Proposition}
\newtheorem{remark}[theorem]{Remark}
\newtheorem{definition}[theorem]{Definition}

 \def\bbbz{{\mathchoice {\hbox{$\sf\textstyle Z\kern-0.4em Z$}}
 {\hbox{$\sf\textstyle Z\kern-0.4em Z$}}
 {\hbox{$\sf\scriptstyle Z\kern-0.3em Z$}}
 {\hbox{$\sf\scriptscriptstyle Z\kern-0.2em Z$}}}}

 \def\bbbc{{\mathchoice {\setbox0=\hbox{$\displaystyle\rm C$}\hbox{\hbox
 to0pt{\kern0.4\wd0\vrule height0.9\ht0\hss}\box0}}
 {\setbox0=\hbox{$\textstyle\rm C$}\hbox{\hbox
 to0pt{\kern0.4\wd0\vrule height0.9\ht0\hss}\box0}}
 {\setbox0=\hbox{$\scriptstyle\rm C$}\hbox{\hbox
 to0pt{\kern0.4\wd0\vrule height0.9\ht0\hss}\box0}}
 {\setbox0=\hbox{$\scriptscriptstyle\rm C$}\hbox{\hbox
 to0pt{\kern0.4\wd0\vrule height0.9\ht0\hss}\box0}}}} 

\newcommand{\nc}{\newcommand}
\nc{\cH}{{\mathcal H}}
\nc{\cA}{{\mathcal A}}
\nc{\cG}{{\mathcal G}}
\nc{\cC}{{\mathcal C}}
\nc{\cO}{{\mathcal O}}
\nc{\cI}{{\mathcal I}}
\nc{\cB}{{\mathcal B}}
\nc{\cY}{{\mathcal Y}}
\nc{\cK}{{\mathcal K}}
\nc{\cX}{{\mathcal X}}
\nc{\cS}{{\mathcal S}}
\nc{\cE}{{\mathcal E}}
\nc{\cF}{{\mathcal F}}
\nc{\cZ}{{\mathcal Z}}
\nc{\cQ}{{\mathcal Q}}
\nc{\cN}{{\mathcal N}}
\nc{\cP}{{\mathcal P}}
\nc{\cL}{{\mathcal L}}
\nc{\cM}{{\mathcal M}}
\nc{\cT}{{\mathcal T}}
\nc{\cW}{{\mathcal W}}
\nc{\cU}{{\mathcal U}}
\nc{\cD}{{\mathcal D}}
\nc{\cJ}{{\mathcal J}}
\nc{\cV}{{\mathcal V}}
\nc{\fr}{{\rightarrow}}

\begin{document}

\title{Surfaces with $p_{g}=q=2$ and an irrational pencil}
\author{Francesco Zucconi\footnote{Mathematical subject classification 14J29, 
14J25, 14D06, 14D99. Partially supported by 
SC.D.I.M.I. cecu 04118 '99 Ricerca Dipartimentale,
Universit\`{a} di Udine.}}
\date{}
\maketitle
\thispagestyle{empty}

\begin{abstract}
We classify all the irrational pencils over the surfaces of general type 
with $p_{g}=q=2$; as a byproduct 
it gives an evidence 
for the Catanese conjecture on surfaces of general type 
with $p_{g}=q=2$. 
\end{abstract}      
{\begin{center}{{\bf{INTRODUCTION}}}\end{center}}

Catanese conjecture about surfaces $X$ of general type 
with $p_{g}=q=2$ states that if $X$ has no irrational 
pencil then $X$ is the double cover of a principally 
polarized Abelian surface branched on a divisor $D$ 
linearly equivalent to $2\Theta$. 
This conjecture is the analogue in the case 
$p_{g}=q=2$ of a similar and known property in the case $p_{g}=q=3.$ 

In fact in [CCM] a partial classification was given of surfaces with
$p_g = q=3$ and recently in \cite{Pi}, \cite{HP} it is showed that 
there are only two families of surfaces, described
in [CCM], with those invariants. For 
related problems on surfaces 
with $p_{g}=q=2$ we refer to \cite{Ci} and 
\cite{CM}.

In this paper we are concerned 
on the other side of the theory. 
In fact we classify all
irrational pencils over 
surfaces of general type with $p_{g}=q=2$. 
The basic result is the following:\\\\
\noindent
{\bf{Theorem [1]}}
{\it{Every irrational pencil of curves with genus $g>2$ 
over a surface of Albanese general type is isotrivial if 
$p_{g}=q=2.$}}\vspace{.10in}\vspace{.10in}

The 
technical point relies 
on the famous result by Simpson \cite{Sim} 
about the locus $V^{1}(X)=\{L\in{\rm{Pic}}^{0}(X)\mid h^{1}(X, -L)\geq 0\}$. 
It follows, almost verbatim, like the  
application of this important technique to 
the surfaces theory which is in \cite{Z1}; incidentally \cite{HP} has a 
similar approach.
Then we classify irrational isotrivial pencil of genus $g>2$ and we obtain:\\

\noindent
{\bf{Theorem [2]}}
{\it{ $S$ is a surfaces of Albanese general type with 
$p_{g}=q=2$ having an irrational 
isotrivial pencil of genus $g>2$ if and only if 
it is the minimal desingularization of 
$C_{1}\times C_{2}/G$ where $G=\mathbb Z/2 \times \mathbb Z/2$, 
$C_{1}$, $C_{2}$ 
have genus $3$, $C_{i}\rightarrow E_{i}=C_{i}/G$ is a 
Galois covering branched over two points of the elliptic curve 
$E_{i}$ and $G$ acts diagonally. 
In this case the general fiber $F$ of the canonical morphism is obtained 
by the smoothing of two curves of genus $3$ which intersect in four  points, 
then it has genus $9$.}}\vspace{.10in}\vspace{.10in}

This theorem needs the 
technique built in \cite{Ca1} and \cite{Z2} 
to deal with isotrivial fibrations. 
In the genus-$2$ case, the description of $S$ as a double covering 
of $Y=E_{1}\times E_{2}$  
branched over 
$D\in \mid 2\pi_{1}^{\star}(P_{1})+ 2\pi_{2}^{\star}(P_{2})\mid$, 
where $E_{i}$, $i=1,2$ are two elliptic curves 
$\pi_{i}:Y\rightarrow E_{i}$ are the natural projections, 
$P_{i}\in E_{i}$ and 
$D$ is a reduced divisor, is not hard to obtain. 
Here we recall only the result:\\\\
\noindent
{\bf{Proposition [3]}}
{\it{$S$ is a surfaces of Albanese general type with 
$p_{g}=q=2$ having an irrational pencil of genus $2$
if and only if $S$ is 
the normalization of the double cover of $Y$ branched over $D$.
In particular the general fiber $F$ of the canonical map 
$\phi_{\mid K_{S}\mid}:S-\rightarrow \mathbb 
P^{1}$ is a genus $5$ curve.}}
\vspace{.10in}\vspace{.10in}

\noindent
A glance to our results show that we have an 
indirect evidence for the Catanese conjecture. 
The non isotrivial fibrations appear only in 
the degenerate case of the conjecture. 
A way to prove the conjecture could be to produce 
non isotrivial fibrations on $S$ by a degeneration argument.

The last step is 
the classification of surfaces $S$ with $p=q=2$ not of Albanese general type.
In this case the Albanese morphism induces a fibration 
$\psi:S\rightarrow C\subset{\rm{Alb}}(S)$ over a curve of genus $2$ 
and $\psi$ is an \'{e}tale bundle of genus $2$-fiber $F$. 
Then $S$ is the quotient of $C^{'}\times F$ by the diagonal $G$-actions where 
$F/G=\mathbb P^{1}$ and $C^{'}\rightarrow C=C^{'}/G$ is \'{e}tale. 
These surfaces are the first occurrences of the 
concept of G(eneralized)H(yperelliptic)- 
Surface cf.\cite{Ca1}, \cite{Z3}: see proposition \ref{ipgeneralizzate}. 
Nowadays, the irreducible components of the moduli space of GH-surfaces 
have been described, \cite[Theorem B]{Ca1}.   
Then the Bolza classification in \cite{Bl} of $G$-actions on a genus-$2$ 
curve enables to end this classification, see \ref{finiamola}, 
and, simultaneously, it completes the wanted description of 
all the irrational pencils.

It is a pleasure to take this 
opportunity to thank 
Fabrizio Catanese who helped me 
to write a modern proof of the Bolza's classification.
I also thank Ciro Ciliberto for 
the encouragement to write this article and 
Gian Pietro Pirola for 
an enlightening mathematical conversation 
on this topic. Finally, 
I thank my colleagues of D.I.M.I 
in Udine for the extension 
of my teaching-duty free period and 
the {\it{Departamento de 
Matem\`{a}tica Aplicada I de la 
Universitat 
Polit\`{e}cnica de Catalunya}}, where I wrote this article, 
for the invitation 
to stay there during the Autumn 2001.

\section{Poof of theorem[1]}
In this section $S$ is a surface of general 
type with $p_{g}=q=2$ whose Albanese morphism 
$\alpha:S\longmapsto {\rm{Alb}}(S)$ is surjective.

\begin{lemma}\label{cuore}
Let $S$ be a surface of Albanese general type with 
$p_{g}=q=2$. Let $\phi:S\rightarrow B$ be a fibration of genus $g$. 
If the genus $b$ of $B$ is $>0$ then $b=1$ and 
$$
\phi_{\star}\omega_{S}=
{\cal {O}}_{B}\bigoplus_{i=1}^{g-2} {\cal {O}}_{B}(\eta_{i})\bigoplus \cL
$$
\noindent
where $\eta_{i}\in{\rm{Tors}}({\rm{Pic}}^{0}(B))\setminus\{0\}$ and $\cL$ 
is an invertible sheaf of degree $1$.
\end{lemma}
\proof If $b=2$ then $S$ is not of Albanese general type by 
the universal property of $\alpha$. Then $B$ is an elliptic curve. 
By \cite{Fu}, $\phi_{\star}\omega_{S}=
{\cal {O}}_{B}\bigoplus{\cal F}$ where $\cal F$ is a 
nef locally free
sheaf of rank $g-1$ 
such that $h^{0}(B,\, {\cal F}^{\vee})=0$. Let 
${\cal F}=\bigoplus_{i=1}^{k}{\cal F}_{i}$ be the decomposition into 
indecomposable subvectorbundles. Obviously $h^{0}(B,{\cal F})=1$. Then 
by Riemann-Roch on $B$, ${\rm{deg}} {\cal F}=1$.

Let ${\cal L}$ be the smallest subvectorbundle of ${\cal F}$ containing 
the subsheaf generically generated by $H^{0}(B,{\cal F})$. 
In particular ${\rm{deg}} {\cal L}>0$ and it has rank $1$; that is,  
$$
{\cal L}\sim{\cal {O}}_{B}(P)
$$\noindent where $\sim$ means linear equivalence.

\noindent{\it{Claim:}} {\it{up to reorder the ${\cal F}_{i}$'s, it holds 
${\cal L}={\cal F}_{1}$.}}

$\subset .$ The inclusion 
$\phi:{\cal L}\hookrightarrow{\cal F}$ gives for each factor the morphism 
$\pi_{i}\circ\phi= \phi_{i}:{\cal L}\rightarrow{\cal F}_{i}$. Set 
$\phi_{i}({\cal L})= {\cal L}_{i}$ and let $d_{i}= {\rm{deg}} {\cal L}_{i}$. 
If ${\cal L}_{i}\neq 0$ then $d_{i}\geq 1$. Since $1=h^{0}(B,{\cal F})\geq 
\sum_{i=1}^{k}d_{i}$ then, up to reorder, we have $d_{1}=1$, $d_{i}=0$ and 
${\cal L}_{i}=0$ where $i=2,...,k$. 

$\supset .$ By contradiction. 
Assume that ${\cal F}_{1}\neq {\cal L}$. By definition of ${\cal L}$ the 
quotient ${\cal F}_{1}/{\cal L}$ is locally free. Let 

$$
{\cal F}_{1}/{\cal L}=\bigoplus_{i=1}^{m}{\cal F}_{1,i}
$$\noindent
be a direct sum of its indecomposable components. 
Note that: ${\rm{deg}} {\cal F}_{1,i}=0$. Then 

\begin{equation}\label{barrow}
H^{1}(B,{\cal L}\otimes(\bigoplus_{i=1}^{m}{\cal F}_{1,i})^{\vee})=
\bigoplus_{i=1}^{m}
H^{1}(B,{\cal L}\otimes({\cal F}_{1,i})^{\vee})=
\bigoplus_{i=1}^{m}
H^{0}(B,{\cal L}^{\vee}\otimes {\cal F}_{1,i})=0.
\end{equation}
Now (\ref{barrow}) and \cite[III.6.2]{Ha} imply that 
$$
0\rightarrow {\cal L}
\rightarrow {\cal F}_{1}
\rightarrow\bigoplus_{i=1}^{m}{\cal F}_{1,i}  \rightarrow 0
$$\noindent splits: a contradiction, since ${\cal F}_{1} $ is indecomposable.

\noindent{\it{Claim:}} 
{\it{${\cal F}_{j}\in{\rm{Tors}}({\rm{Pic}}^{0}(B))\setminus\{0\}$ 
if $j=2,...,k.$}}

By Atiyah classification of vector bundles over an elliptic curve, \cite{At}, 
we have ${\cal F}_{j}=F_{r_{j}}\otimes T_{j}$ where the rank-$r_{j}$ sheaves 
$F_{r_{j}}$ are obtained inductively through non trivial extensions by 
$\cO_{B}$. Besides $T_{j}\in{\rm{Pic}}^{0}(B)$. 
We want to prove that $T_{j}$ is torsion and $r_{j}=1$.

\noindent
{\it{First step: $T_{j}$ is torsion.}}

Let $\Lambda\in {\rm{Pic}}^{0}(B)\setminus\{0\}.$ By the 
Serre duality and by the projection formula it holds:

\begin{equation}\label{tesisi}
h^{2}(S, K_{S}+f^{\star}(\Lambda))=
h^{0}(S, -f^{\star}(\Lambda))=h^{0}(E, -\Lambda)=0.
\end{equation}\noindent
By the Leray spectral sequence for the morphism $f$, 
$$h^{1}(S, K_{S}+f^{\star}(\Lambda))=h^{1}(B,\Lambda)\oplus 
h^{1}(B, \cL\otimes \Lambda)\oplus_{j=1}^{n}h^{1}(B, 
{\cal F}_{j}\otimes \Lambda)
+h^{0}(B,\Lambda).
$$\noindent 
Then by Riemann-Roch on 
$B$ and by relative duality
we have:

\begin{equation}\label{tesisini}
h^{1}(S, K_{S}+f^{\star}(\Lambda))=
\oplus_{j=1}^{n}h^{1}(B, {\cal F}_{j}\otimes \Lambda).
\end{equation}\noindent
Choose $\Lambda=-T_{i}$ where $2\leq i\leq k$. 
Then we have a jump in cohomology and by the Simpson solution of the Beauville-Catanese conjecture 
\cite{Sim}, it implies that $T_{i}$ are torsion sheaves.

\noindent
{\it{Second Step: $r_{j}=1$.}}

Only for clarity reasons we assume that the torsion sheaves 
are of relative prime order. By contradiction assume that there exists 
$2\leq i\leq k$ with $r_{i}>1$.

Let $\tau_{i}:E_{i}\rightarrow B$ 
be the unramified covering given by $T_{i}$ and denote by
$f_{i}:S_{i}=S\otimes_{B}E_{i}\rightarrow E_{i}$ and by 
$\sigma_{i}:S_{i}\rightarrow S$ the two projections. Then 
$\tau_{i}^{\star}{\cal F}_{i}=F_{r_{i}}$, 
$\omega_{S_{i}}=\sigma_{i}^{\star}\omega_{S}$ and  
$f_{i\star}\sigma^{\star}_{i}\omega_{S}=
\tau_{i}^{\star}f_{\star}\omega_{S}.$ 
In particular, $F_{r_{i}}$ is a direct summand of 
$f_{i\star}\omega_{S_{i}}$. Then 
$h^{1}(E_{i}, f_{i\star}\omega_{S_{i}})=2$. By 
\cite{Fu}, $\cO_{B}$ must be a direct summand of 
$f_{i\star}\omega_{S_{i}}$; this forces 
$F_{r_{i}}$ to be decomposable: a contradiction.\qed\vspace{.10in} 

\begin{corollary}\label{irratio}
Any irrational pencil 
over 
$S$ is an elliptic 
pencil of genus 
$2\leq g\leq 5$.
\end{corollary}
\proof By Stein factorization, every irrational 
pencil gives an elliptic fibration; choose $\phi:S\rightarrow B$ 
one of these fibrations.
Following \ref{cuore} even in the notation, 
we have a direct summand 
${\cal L}\hookrightarrow\phi_{\star}\omega_{S}$ 
such that, 
${\cal L}\sim{\cal {O}}_{E}(P)$, $P\in B$. 
Then by Xiao's method cf.\cite{BZ}, it follows that $K_{S}-D$ is nef 
where $D=\phi^{-1}(P)$. Then $(K_{S}-D)K_{S}\geq 0$. In particular by 
Miyaoka's inequality $9\geq K_{S}^{2}\geq K_{S}D$. Since $b>0$ then $D^{2}=0$ and the inequality $2\leq g\leq 5$ now easily follows 
by the adjunction formula.\qed\vspace{.10in}

Combining the universal property of the Albanese morphism and a standard 
Abelian varieties argument, once we have an elliptic pencil 
$f_{1}:S\rightarrow E_{1}$ 
then there exists another elliptic pencil (not necessarily distinct) 
$f_{2}:S\rightarrow E_{2}$.
We recall also that $\xi:A^{'}\rightarrow A$ is a {\it{minimal isogeny}} 
among two Abelian varieties if every isogeny 
$A^{'}\rightarrow A^{''}$ factorizes through $\xi$.

\begin{corollary}\label{appena}
There exists a minimal isogeny 
$j$ of ${\rm{Alb(S)}}$ over a product of two elliptic curves
$E_{1}\times E_{2}$.  
\end{corollary}
\proof Trivial.\qed\vspace{.10in}

\begin{definition}\label{perpoi}
The two elliptic fibration $f_{i}:S\rightarrow E_{i}$, $i=1,2$ 
induced by $j$ will be called natural fibrations of 
$S$ and $g_{i}$ will denote the genus of the general $f_{i}$-fiber.  
\end{definition}
For further reference we sum up some results on irrational pencils:

\begin{proposition}\label{finale}
If a surface of Albanese general 
type with $p_{g}=q=2$ 
has an irrational pencil then ${\rm{Alb(S)}}$ 
is minimally 
isogenus to 
a product of two elliptic curves 
$E_{1}\times E_{2}$. 
Every irrational pencil 
factorizes through 
one of the natural 
fibrations $f_{i}$, $i=1,2$. 
Moreover it holds $2g_{i}-2\geq F_{1}F_{2}$.
\end{proposition}
\proof We have to show only the inequality. 
By the proof of \ref{irratio} we know that $K_{S}-F_{j}$ is nef. Then 
$2g_{i}-2=K_{S}F_{i}\geq F_{j}F_{i}$, $i,j\geq 1,2$. \qed\vspace{.10in}

The following theorem actually shows that all the Jacobian of the
irreducible components of an irrational pencil dominates a 
``large'' Abelian variety; this forces isotriviality.  

\begin{theorem}\label{deformazione}
Every elliptic fibration with fiber of genus $g>2$ over a surface 
with $p_{g}=q=2$ is an isotrivial pencil.
\end{theorem}
\proof By the universal property of the Albanese morphism such 
surface is of Albanese general type.
Let $f_{1}:S\rightarrow E_{1}$ be one of the two natural 
fibrations of definition \ref{perpoi} and assume that $g_{1}>2$. 
By \ref{cuore} it holds:
\begin{equation}\label{barcellona}
f_{1\,\star}\omega_{S}=
{\cal {O}}_{E_{1}}\bigoplus_{i=1}^{g-2} {\cal {O}}_{E_{1}}
(\eta_{i})\bigoplus \cL
\end{equation}
\noindent
where $\eta_{i}\in{\rm{Tors}}({\rm{Pic}}^{0}(E_{1}))\setminus{0}$ and $\cL$ 
is an invertible sheaf of degree $1$. Let $\sigma:E\rightarrow E_{1}$ 
be the unramified base change given 
by ${\rm{lcm}}\{\eta_{i}\}_{i=1}^{g-2}$. 
$X=S\otimes_{E_{1}}E$ is connected. Let $\psi:X\rightarrow E$, 
$\tau:X\rightarrow S$ be the projections. It is easy to 
see that from (\ref{barcellona}):
\begin{equation}\label{cellona}
\psi_{\star}\omega_{E}=
\bigoplus_{i=1}^{g-1} {\cal {O}}_{E}\bigoplus \tau^{\star}\cL.
\end{equation}\noindent By \cite{Fu} it follows that $q(X)=g$. 
By the universal property, the Jacobian over 
the smooth fiber has a surjection 
$\mu_{t}:J(F_{t})\rightarrow A$ where $t\in E$, $F_{t}$ is smooth and $A$ 
is a {\it{fixed}} Abelian subvariety of $\rm{Alb}(X)$. 
In particular, up to isogenies over $A$, there exists an elliptic curve 
$E_{t}\hookrightarrow J(F_{t})$ such that 
\begin{equation}\label{barc}
0\rightarrow E_{t}\rightarrow J(F_{t})\rightarrow A\rightarrow 0
\end{equation}\noindent is exact. Now if $E_{t}$ is independent of 
$t$ then the Abelian variety rigidity imply that $J(F_{t})$ is fixed. 
Then the smooth fibers are isomorphic. This means that $f_{1}$ is isotrivial.

Assume that the $E_{t}$'s are a non constant family. By \ref{finale} 
the generic smooth $F_{t}$ is equipped with two morphisms: 
$F_{t}\rightarrow E_{2}$, $F_{t}\rightarrow E_{t}$ 
where $E_{2}$ is the fixed elliptic curve which is 
the basis of the other natural fibration 
$f_{2}:S\rightarrow E_{2}$.  
Since $E_{t}$ moves in a continuous family, the induced 
morphism into the product 
$\nu_{t}:F_{t}\rightarrow E_{2} \times E_{t}$ is an embedding. 
By adjunction, it easily follows $g(F_{t})\geq 5$. 
Then $g(F_{t})=5$, since \ref{irratio}. In particular 
$F_{t}\rightarrow E_{2}$ is a $2$-to-$1$ morphism. 
The same argument works interchanging the role between $f_{1}$ and $f_{2}$. 
Then $S\rightarrow E_{1}\times E_{2}$ is a generically finite 
$2$-to-$1$ covering whose branch locus $\Delta\in\mid 2\delta\mid$ 
satisfies the conditions $\delta \, (E_{1}\times\{y\})=4$,
$\delta \, (\{x\}\times E_{2}) =4$ where $x\in E_{1}$ and  
$y\in E_{2}$. On the other hand $\delta^{2}\leq 2$ since 
$9\geq K_{S}^{2}\geq 4\delta^{2}$. It is easy to see that 
$\delta$ does not exist.\qed\vspace{.10in} 

We have shown the theorem[1]. For further use we prove:
\begin{proposition}\label{diagonale}
Every surface $S$ of Albanese general type with $p_{g}=q=2$ 
equipped with an irrational pencil with fiber of genus $g>2$ 
can be realized 
as the  minimal desingularization 
of the quotient surface $X/G$, where 
$X=B_{1}\times B_{2}$, $G$ acts 
diagonally over $X$ and 
faithfully over the smooth curves $B_{1}$, $B_{2}$. 
Moreover $B_{1}/G=E_{1}$, $B_{2}/G=E_{2}$ where 
$E_{1}$, $E_{2}$ are the basis of the two normal fibrations 
over $S$.
\end{proposition}
\proof
By Stein factorization the pencil induces on 
$S$ an elliptic fibration. By \ref{deformazione} and 
by \cite{Ser} the claim follows.\qed\vspace{.10in}

\section{Proof of theorem[2]}
By \ref{diagonale} the 
classification task is reduced to classify all the Galoisian $G$ 
actions over curves $B_{1}$ and $B_{2}$ of genus at most $2\leq b_{i}\leq 5$, 
$i=1,2$ such that the quotient curves $B_{1}/G=E_{1}$, $B_{2}/G=E_{2}$ 
are elliptic and the 
diagonal $G$ action over $X=B_{1}\times B_{2}$ has quotient with $p_{g}=q=2$. 
If $G$ is {\it{Abelian}} let $G^{\star}$ be the group of the characters.

\begin{lemma}\label{allafine}
Every surface $S$ of Albanese general type with $p_{g}=q=2$ 
equipped with an isotrivial pencil can be realized 
as the  minimal desingularization of the quotient surface $X/G$ 
where $G$ is Abelian
if and only if the faithful $G$ action on the two smooth curves 
$B_{1}$, $B_{2}$ 
of genus $2\leq b_{i}\leq 5$, $i=1,2$ gives two decompositions 
$H^{0}(B_{1}, \Omega^{1}_{B_{1}})=
\bigoplus_{\chi\in G^{\star}}V^{(1)}_{\chi}$, 
$H^{0}(B_{2}, \Omega^{1}_{B_{2}})=
\bigoplus_{\chi\in G^{\star}}V^{(2)}_{\chi}$ 
such that there exists a unique 
nontrivial character $\chi$ with 
$V^{(1)}_{\chi}\otimes V^{(2)}_{\chi^{-1}}\neq 0$. Moreover 
the following two numerical conditions hold: 
$1)\,\, {\rm{dim}}_{\mathbb C}V^{(1)}_{{\rm{id}}}= 
{\rm{dim}}_{\mathbb C}V^{(2)}_{{\rm{id}}}=1$ and 
$2)\,\, {\rm{dim}}_{\mathbb C}V^{(1)}_{\chi}=
{\rm{dim}}_{\mathbb C}V^{(2)}_{\chi^{-1}}=1.
$
\end{lemma}
\proof The direct proof is easy; otherwise 
cf.\cite[Lemma 1.3]{Z2} and \cite[Theorem 1.4]{Z2}.\qed\vspace{.10in}

We need to compute the $G$-actions over $F$ where 
$2\leq g(F)\leq 5$ and $F/G$ is elliptic.

\begin{lemma}\label{azioni}
Let $\pi:F\rightarrow E$ be a Galois morphism with group $G$ 
such that $E$ is an elliptic curve and $F$ has genus $2\leq g\leq 5$. 
Then the occurring actions are the following:
$$\begin{tabular}{|l|c|l|}\hline
$ g$&$ G$ &$\bigoplus_{\chi\in G^{\star}\setminus{\rm{id}}} V^{i}_{\chi}$ 
\\ \hline
$2$&$ \mathbb Z/2$&
$ V^{1}_{-}$ \\\hline
$3$&$ \mathbb Z/2$&
$V^{2}_{-}$ \\\hline
$3$&$ \mathbb Z/3$&
$V^{1}_{\chi}\oplus V^{1}_{\chi^{2}}$ \\\hline
$3$&$ \mathbb Z/4$&
$V^{1}_{\chi}\oplus V^{1}_{\chi^{3}}$ \\\hline
$3$&$ \mathbb Z/2\times \mathbb Z/2 $&
$V^{1}_{\chi_{1}}\oplus V^{1}_{\chi_{12}}$ \\\hline
$4$&$ \mathbb Z/2 $&
$V^{3}_{-}$ \\\hline
$4$&$ \mathbb Z/3 $&
$V^{1}_{\chi}\oplus V^{2}_{\chi^{2}}$\\\hline
$4$&$ \mathbb Z/2\times\mathbb Z/2  $&
$V^{1}_{\chi_{1}}\oplus V^{1}_{\chi_{12}}\oplus V^{1}_{\chi_{2}}$\\\hline
$4$&$ \mathbb Z/4$&
$V^{1}_{\chi}\oplus V^{1}_{\chi^{2}}\oplus V^{1}_{\chi^{3}}$\\\hline
$4$&$ \mathbb Z/6$&
$V^{1}_{\chi}\oplus V^{1}_{\chi^{3}}\oplus V^{1}_{\chi^{5}}$\\\hline
$4$&$ {{\cS}}_{3}$&$U^{2}\oplus W^{1}$\\\hline
$4$&$ {{\cS}}_{3}$& $W^{3}$\\\hline
$5$&$ \mathbb Z/2 $&
$V^{4}_{-}$ \\\hline
$5$&$ \mathbb Z/3 $&
$V^{2}_{\chi}\oplus V^{2}_{\chi^{2}}$\\\hline
$5$&$ \mathbb Z/2\times\mathbb Z/2  $&$
V^{2}_{\chi_{1}}\oplus V^{2}_{\chi_{12}}$\\\hline
$5$&$ \mathbb Z/2\times\mathbb Z/2  $&$
V^{2}_{\chi_{1}}\oplus V^{1}_{\chi_{12}}\oplus V^{1}_{\chi_{2}}$ \\\hline
$5$&$ \mathbb Z/4$& $V^{2}_{\chi} \oplus 
V^{1}_{\chi^{2}}\oplus V^{1}_{\chi^{3}}$\\\hline
$5$&$ \mathbb Z/4$& $V^{2}_{\chi} \oplus V^{2}_{\chi^{3}}$\\\hline
$5$&$ \mathbb Z/5$&
$\oplus_{i=1}^{4}V^{1}_{\chi^{i}}$\\\hline
$5$&$ {{\cS}}_{3}$&$U_{1}^{2}\oplus U_{2}^{2}$\\\hline
$5$&$ {{\cS}}_{3}$&$U^{2}\oplus W^{2}$\\\hline
$5$&$ {{\cS}}_{3}$& $W^{4}$\\\hline
$5$&$ (\mathbb Z/2)^{3} $&
$V^{1}_{\chi_{1}}\oplus V^{1}_{\chi_{2}}\oplus V^{1}_{\chi_{3}}\oplus 
V^{1}_{\chi_{123}}$ \\\hline
$5$&$ \mathbb Z/2\times \mathbb Z/4 $&
$V^{1}_{\chi_{1}}\oplus V^{1}_{\chi_{1}\chi_{2}}
\oplus V^{1}_{\chi_{1}\chi_{2}^{2}}\oplus V^{1}_{\chi_{1}\chi_{2}^{3}}$\\\hline
$5$&$ \mathbb Z/2\times \mathbb Z/4 $&
$V^{1}_{\chi_{1}}\oplus V^{1}_{\chi_{2}}
\oplus V^{1}_{\chi_{2}^{3}}\oplus V^{1}_{\chi_{1}\chi_{2}^{2}}$\\\hline
$5$&$ \mathbb Z/2\times \mathbb Z/4 $&
$V^{1}_{\chi_{2}}\oplus V^{1}_{\chi_{2}^{3}}
\oplus V^{1}_{\chi_{2}}\oplus V^{1}_{\chi_{1}\chi_{2}^{2}}$\\\hline
$5$&$ \mathbb Z/8$&
$V^{1}_{\chi}\oplus V^{1}_{\chi^{3}}
\oplus V^{1}_{\chi^{5}}\oplus V^{1}_{\chi^{7}}$\\\hline
$5$&$ \mathbb \cD_{4}$&
$V^{2}_{\chi_{1}}\oplus V^{2}_{\chi^{12}}$\\\hline
\end{tabular}$$\noindent
where $V^{i}_{\chi}$ means that the $\chi$-piece of 
$H^{0}(F, \Omega^{1}_{F})$ has dimension $i$. For the 
${{\cS}}_{3}$-representations we have denoted by $U^{2}$, $U^{2}_{i}$ 
the irreducible subspaces of dimension $2$ and for the 
$\cD_{4}$-representations only (two of) the {\underline{linear}} 
characters occur.
\end{lemma}
\proof It is an application of the Riemann Hurwitz formula plus a 
careful analysis on the action over the branch loci. 
Here we show how the quaternion actions over a genus $5$ 
curve with elliptic quotient can be excluded. 
In fact the orbifold exact sequence 
cf.\cite[definition 4.4]{Ca1} for these actions is:

$$0\rightarrow 
\pi(F)\rightarrow\langle a,b,x,y\mid x^{2}=y^{2}=xy[ab]=1\rangle 
\stackrel{\mu}{\rightarrow}\cH\rightarrow 0
$$\noindent 
where $\cH$ is the quaternion group. Since 
$\cH$ has only one element of order $2$, denote it by $-1,$ 
then $\mu: x\mapsto -1$ and $\mu: y\mapsto -1$. In particular 
$[\mu (a),\mu(b)]=1$ and $\mu$ cannot be surjective: a contradiction. 
Using the fact that the dihedral group of order $8$, $\cD_{4}$, 
has one normal subgroup of order $2$ we can compute also this action.\qed
\vspace{.10in}

Now the final step is to use \ref{diagonale} and 
\ref{allafine} to construct the quotient surface. 

{\begin{center}{\bf{Construction}}\end{center}}

To exclude the non Abelian cases it requires only to 
couple the possible actions to see that it never happens that $p_{g}=q=2$. 
The computation for the Abelian case can be easily done.
The two solutions correspond to a  ${\mathbb Z/2}$ diagonal 
action on the product of two genus-$2$ curves and to a diagonal
$G=\mathbb Z/2  \times \mathbb Z/2$ action over the 
product of two genus $3$ curves; in this last case 
$G$ acts on the two factors 
via the action $V^{1}_{\chi_{1}}\oplus V^{1}_{\chi_{2}}$ 
and respectively $V^{1}_{\chi_{1}}\oplus V^{1}_{\chi_{12}}$. 
A more geometrical construction can be achieved following \cite{Z2}. 
We have:

\begin{theorem}\label{principale}
There are only two classes $\cM_{\mathbb Z/2}$, 
$\cM_{\mathbb Z/2  \times \mathbb Z/2}$ of Albanese general type surfaces 
with $p_{g}=q=2$ and with an isotrivial pencil. An element 
$S\in \cM_{\mathbb Z/2}$ is the minimal desingularization of the quotient 
surface $B_{1}\times B_{2}/G$ where $G=\mathbb Z/2$, 
$B_{1}$, $B_{2}$ have genus  
$b_{1}=b_{2}=2$ and $B_{i}\rightarrow E_{i}=B_{i}/ G$ is a 
Galois covering branched over two points of the elliptic curve $E_{i}$. 
In this case the general fiber $F$ of the canonical morphism is obtained 
by the smoothing of two curves of genus $2$ which intersect 
in two points, then it  has genus $5$. An element $S\in
\cM_{\mathbb Z/2  \times \mathbb Z/2}$ is the minimal desingularization of 
$C_{1}\times C_{2}/G$ where $G=(\mathbb Z/2)^{2}$, $C_{1}$, $C_{2}$ 
have genus  
$c_{1}=c_{2}=3$ and 
$C_{i}\rightarrow E_{i}=C_{i}/ G$ is a 
Galois covering branched over two points of $E_{i}$. 
In this case the general fiber $F$ of the canonical map 
is obtained 
by the smoothing of two curves of genus $3$ 
which intersect in four points, 
then it has genus $9$. 
\end{theorem}\vspace{.10in}\vspace{.10in}\vspace{.10in}\vspace{.10in}

In particular we have proved the theorem[2].


\section{Surfaces with $p_{g}=q=2$ and non surjective Albanese morphism}
In this section $S$ will be a surface of general type 
with $p_{g}=q=2$ 
and $\alpha$ non surjective. We will show that 
$\psi:S\rightarrow C\subset{\rm{Alb}}({S})$ is 
an \'{e}tale bundle and that $S$ is 
a Generalized Hyperelliptic surface. 
In fact these surfaces are the baby examples of $GH$-surfaces.
\vspace{.10in}\vspace{.10in}
{\begin{center}{\bf{Generalized Hyperelliptic Surfaces}}\end{center}} 
The following definition is in \cite{Ca1}, see also \cite{Z3}. 
Let $C_{1}$, $C_{2}$ be two smooth
curves with the corresponding automorphisms groups: 
${\rm{Aut}}(C_{1})$, ${\rm{Aut}}(C_{2})$. Let $G$ be 
a non trivial finite group $G$ with two injections:
$G\hookrightarrow{\rm{Aut}}(C_{1})$, 
$G\hookrightarrow{\rm{Aut}}(C_{2})$.

\begin{definition}\label{iperellittichegeneralizzate}
The quotient surface $S=C_{1}\times C_{2}/G$
by the diagonal action
of $G$ over $C_{1}\times C_{2}$ is said
to be
of {\rm{Generalized
Hyperelliptic type}} (GH, to short) if \\
$i)$ the Galois morphism $\pi_{1}:C_{1}\rightarrow C=C_{1}/G$ is
unramified;\\
$ii)$ the quotient curve $C_{2}/G$ is
isomorphic to $\mathbb P^{1}.$
\end{definition}

\begin{proposition}\label{ipgeneralizzate}
If $S$ is not of Albanese general type with $p_{g}=q=2$ 
then $S$ is GH.
\end{proposition}
\proof
Since  
$\alpha:S\longmapsto {\rm{Alb}}(S)$ is not surjective, 
it induces a fibration 
$\psi:S\rightarrow C\subset{\rm{Alb}}{S}$ 
where $C$ is a smooth curve of genus $2$. 
By \cite{Fu}  
applied to $\psi$, $\chi (S) \geq g-1$. 
Since $\chi (S)=1$ then $g\leq 2$, but ${\rm{Kod}}(S)=2$ implies $g=2$. 
Then the equality 
$\chi (S)=g-1$ holds, and by 
\cite[ch. III prop. 18.2 and prop. 18.3]{BPV}
$\psi$ is an \'{e}tale bundle with fiber of genus $2$. Then 
there exists a group $G$ acting on two curves $C_{1}$, $C_{2}$ 
such that 
$\pi_{1}:C_{1}\rightarrow C_{1}/G=C$ is \'{e}tale, $g(C_{2})=2$ and 
the base change diagram

\begin{equation}\label{etale}
\begin{array}{ccc}
Y &\stackrel{\pi}{\rightarrow} &S\\
\downarrow & &\downarrow\\
C_{1} &\stackrel{\pi_{1}}{\rightarrow} &C\\
\end{array}
\end{equation}
\noindent 
works with $Y=C_{1}\times C_{2}$. In particular 
$g(C_{1})>g(C)=g(C_{2})=2$ and $G$ acts diagonally on $Y$. 
Then $g(C_{2}/G)=0$ and $S$ is GH.\qed\vspace{.10in}

The {\it{Viceversa}} of \ref{ipgeneralizzate} holds in a strong form:

\begin{proposition} \label{quattro}
Let $C_{1}$, $C_{2}$ and $G$ as above. 
If $g(C_{2})=2$, $C_{2}/G=\mathbb P^{1}$ and 
$\pi_{1}:C_{1}\rightarrow C_{1}/G=C$ is an \'{e}tale morphism where 
$g(C)=2$ then the quotient $S=C_{1}\times C_{2}/G$ by the 
diagonal action is a minimal smooth surface of general type 
with $p_{g}(S)=q(S)=2$ and non surjective Albanese morphism.
\end{proposition}
\proof Since $G$ acts freely over $Y$ then $S$ is minimal, smooth, 
of general type and 
$\chi(Y)=n\chi(S)$, $K_{S}^{2}=1/n K_{Y}^{2}$ being $n$ the order of $G$. 
By Riemann-Hurwitz $g(C_{1})=n+1$ then by 
\cite[III.22]{Be}, $p_{g}(S^{'})=2(n+1)$ and $q(Y)=n+3$; that is $\chi(Y)=n$. Then $\chi(S)=1$ and by \cite[VI.12.1]{Be} it follows $p_{g}(S)=q(S)=2$ since 
$g(F/G)=0\: , \: g(C)=2.$ Let us consider the natural fibration 
$S\rightarrow C$. Since 

$$\begin{array}{ccc}
S &\stackrel{\alpha}{\rightarrow} &{\rm{Alb}}(S)\\
\downarrow &\searrow &\downarrow\\
C &\rightarrow &J(C)\\ 
\end{array}$$
\noindent
is commutative, then $\alpha$ is not surjective.\qed\vspace{.10in}

From \ref{quattro} we can extract the following property 
which is a standard feature of GH-surfaces.

\begin{remark} \label{cinque}
Let $S$ be a minimal surface with $p_{g}=q=2$ such that 
$\alpha(S)=C$ is a curve. Then ${\rm{Alb}}(S)=J(C)$.
\end{remark}

Notice that since $S$ is GH, 
$S\rightarrow C_{2}/G=\mathbb P^{1}$ 
{\it{is not the canonical map}}.

By \ref{quattro} 
we have 
the following 
classification 
criteria
\begin{remark}\label{sei}
To classify $S$ 
is equivalent 
$i)$ 
to classify all 
the $G$ actions 
over a 
genus $2$ curve 
$C_{2}$ such that 
$C_{2}/G= \mathbb P^{1}$ 
and $ii)$ 
for each occurrence of 
$G$ in $i)$ 
to classify all the 
\'{e}tale morphisms 
$\pi_{1}:C_{1}\rightarrow C_{1}/G=C$ 
where 
$g(C)=2$.
\end{remark} 
We give a 
modern way 
to solve 
$i)$ and 
so we will 
have a 
classification 
of $S$ following 
\cite[Theorem B]{Ca1}.

The classification of all the couples $(C,\, G)$ where 
$C$ is a curve of genus $2$ and $G$ is a 
subgroup of ${\rm{Aut}}(C)$ 
was obtained by Bolza \cite{Bl}. 
However, probably because he assumed 
it to be trivial (at least for our 
understanding of his proof) 
he did not specify the different 
dihedral actions. 
In every case 
we rewrite it again because 
we adopt a completely 
different approach which 
should be generalizable to 
the other hyperelliptic curves 
and also because we give a 
recipe to compute 
easily each action.

{\begin{center}{\bf{Weighted Projective Spaces}}\end{center}}
Let $C$ be a curve of genus $2$. 
The canonical ring 
${\cal{R}}\bigoplus_{m\geq 0}{\cal{R}}_{m}$ of $C$, 
where 
${\cal{R}}_{m}=H^{0}(C,\, \omega _{C}^{\otimes m})$ 
manifests 
$C$ as an hypersurface $C={\rm{Proj}}({\cal{R}})$ of 
the weighted projective space 
$\mathbb P(1,1,3)={\rm{Proj}}\mathbb C[x_{0},x_{1},z]$ 
where $x_{0},x_{1}$ 
have degree $1$ and 
$z$ has degree $3$. 
We want to find explicitly 
the action of ${\rm{Aut}}(C)$ over 
${\rm{Proj}} ({\cal{R}})$. 
This description requires 
some well known facts 
that we recall since 
they can be useful 
for further generalizations.  

The genus $2$ curve $C$ 
is the normalization 
of the projective closure 
${\overline{C}}_{0}\subset \mathbb P^{2}$ of  
$C_{0}=\{(x,y)\in \mathbb A^{2}\mid\, y^{2}=\beta (x)\}$ where 
$\beta\in\mathbb C [x]$, 
${\rm{deg}}\beta=6$ and it has 
$6$ distinct roots. 
Let 
$\nu:C\rightarrow{\overline{C}}_{0}$ 
be the normalization 
morphism then the 
hyperelliptic involution 
$i:C\rightarrow C$ is induced by the affine automorphism
$(x,y)\rightarrow (x,-y)$ and
\begin{equation}\label{forme}
\omega _{0}=
\nu ^{\star}\frac{dx}{y}\,\,\,,\,\,\omega _{1}=\nu ^{\star}x\frac{dx}{y}
\end{equation}
\noindent
give a basis of $H^{0}(C,\, \omega _{C})$. Set 
$\eta=\nu ^{\star}y(\frac{dx}{y})^{3}.$ Then

\begin{equation}\label{formedue}
\langle \, \omega _{0}^{3},\, \omega _{0}^{2}\omega _{1},\,
\omega _{0}\omega _{1}^{2},\,\omega _{1}^{3},\, \eta\, \rangle
\end{equation}
\noindent
is a basis of
$H^{0}(C,\, \omega _{C}^{3}).$ By (\ref{forme}) it holds 
$$i^{\star}\omega _{0}=i^{\star}\circ \nu ^{\star}\frac{dx}{y}=
\nu ^{\star}\circ i^{\star}\frac{dx}{y}=-\nu ^{\star}\frac{dx}{y}
=-\omega _{0}$$ \noindent
and in the same way $i^{\star}\omega _{1}=-\omega _{1}.$ In particular if
$\rho:{\rm{Aut}}(C)\rightarrow{\rm{GL}}( H^{0}(C,\, \omega _{C}))$ is the natural faithful representation we have that $\rho (i)=-Id$ and $i$ 
commutes with every $g\in G$. By the tricanonical morphism there is also a faithful representation 
$\rho_{3}:{\rm{Aut}}(C) \rightarrow GL( H^{0}(C,\, \omega _{C}^{\otimes3}))$.

\begin{lemma}\label{cesare}
The hyperelliptic involution $i$ is 
in the center of ${\rm{Aut}}(C)$. The action of $i$ splits: 
$$H^{0}(C,\, (\omega _{C})^{\otimes 3})={\cal S}^{3}
H^{0}(C,\, \omega _{C})\bigoplus \eta \mathbb C.$$\noindent
Moreover the decomposition is preserved by every $g\in {\rm{Aut}}(C)$.
\end{lemma}
\proof We have just seen that $i$ is central. Obviously 
$i^{\star}(\omega _{0}^{s}\omega _{1}^{j})=-\omega _{0}^{s}\omega _{1}^{j}$ 
with $s+j=3$ and $i^{\star}\eta =\eta;$ and then we have the 
claimed decomposition for $i$. Let $g\in {\rm{Aut}}(C)$ we want to show that 
$g^{\star}\eta =\chi (g)\eta,$ 
where $\chi:{\rm{Aut}}(C)\rightarrow \mathbb C$ is a 
character of ${\rm{Aut}}(C)$. By (\ref{formedue}), 
$g^{\star}\eta=\sum_{i+j=3}a_{ij}\omega _{0}^{i}\omega _{1}^{j}+
\chi (g)\eta.$ Since $g^{\star}\eta=g^{\star}i^{\star}\eta=
i^{\star}g^{\star}\eta$ then $\sum_{i+j=3}a_{ij}\omega _{0}^{i}
\omega _{1}^{j}+\chi (g)\eta=-\sum_{i+j=3}
a_{ij}\omega _{0}^{i}\omega _{1}^{j}+\chi (g)\eta$; that 
is $a_{ij}=0$ for every $i,j$.\qed\vspace{.10in}\noindent

\noindent
Let us consider 
$\mathbb P(1,1,3)={\rm{Proj}}(\mathbb c[x_{0},x_{1},z])$. The map 
$$j: {\cal{R}}\stackrel{}{\longrightarrow}
\frac {{\mathbb  C}[x_{0},\, x_{1},\, z]}{(z^{2}-\beta (x_{0},\, x_{1}))}$$
defined by $\omega _{0}\mapsto x_{0}$, $\omega _{1}\mapsto x_{1}$ 
$\eta \mapsto z$ is an isomorphism; that is 
$C={\rm{Proj}}({\cal{R}})=C_{6}\subset\mathbb P(1,1,3)$. 
In particular, by $j$, we can fix once for 
all an identification $GL(2,\mathbb C)\sim GL(2,\,{\cal R}_{1}).$ 
We easily identify $H^{0}(C,\,\omega_{C})$ and 
$H^{0}(C,\,\omega_{C}^{\otimes 3})$ 
through $C\subset \mathbb P(1,1,3).$

\begin{lemma}\label{augusto}
If $C\subset\mathbb P (1,1,3)$ letting 
$\omega$ the regular differential induced by 
$x_{0}dx_{1}-x_{1}dx_{0}$ it holds: 
\vspace{.10in}\par\noindent
i) $H^{0}(C,\,\omega _{C})=
\{(\frac{\omega}{z})P_{1}(x_{0},x_{1})\mid 
P_{1}(x_{0},x_{1})\in j({\cal{R}})\, ,\,deg(P_{1})=1 \};$
\vspace{.10in}\par\noindent
ii) $H^{0}(C,\,\omega _{C}^{\otimes 3})=
\{ (\frac{\omega}{z})^{3}P_{3}(x_{0},x_{1})\mid P_{3}(x_{0},x_{1})\in 
{\mathbb C}[x_{0},x_{1}] \, ,\,deg(P_{3})=3\,\}.$
\end{lemma}
\proof
A local computation.\qed

\begin{corollary} \label{zurlo}
Via the identification $GL(2,\mathbb C)\sim GL(2,\,{\cal R}_{1})$ it holds:
\vspace{.10in}\par\noindent
$i)$ If $G\subset Aut(C)$ then $G\subset GL(2,\mathbb C);$
\vspace{.10in}\par\noindent
$ii)$ If $\langle i\rangle$ is 
the group generated by the 
hyperelliptic involution then $\langle i\rangle=\{Id,\, -Id\};$
\vspace{.10in}\par\noindent
$iii)$ $G$ acts over $z$ by the character $\chi=det.$
\end{corollary}
\proof Trivial exercise in representation theory.\qed
\begin{corollary} \label{uppa}
$i)$ $\chi(i)=1.$ $ii)$ $G$ acts over $\beta$ via the character 
$\lambda= \chi ^{2}=det ^{2}.$
\end{corollary}
{\begin{center}{\bf{Subgroups of $GL(2,\mathbb C)$}}\end{center}}
We have seen that $G\subset GL(2,\mathbb C).$ Let 
$GL(2, \mathbb C)\stackrel{\pi}{\rightarrow}\mathbb PGL(2, \mathbb C)$ 
be 
the canonical projection,  
$\pi:\, \left( \begin{array}{cc}
a& b\\
c& d
\end{array} \right) \mapsto \left[ \left( \begin{array}{cc}
a& b\\
c& d
\end{array} \right) \right]$ and we set $K=\pi(G).$ We consider the following exact sequences:
\begin{equation}\label{gln}
0\rightarrow \Delta \rightarrow GL(2, \mathbb C)
\stackrel{\pi}{\rightarrow} 
\mathbb PGL(2,\mathbb C)\rightarrow 0
\end{equation}

\begin{equation}\label{pgln}
0\rightarrow \{\pm Id\} \rightarrow SL(2, \mathbb C)
\stackrel{\pi_{1}}{\rightarrow} 
\mathbb PGL(2, \mathbb C)\rightarrow 0.
\end{equation}

\begin{definition}\label{optime}
Let $K$ be a subgroup of 
$\mathbb PGL(2, \mathbb C)$ and let $G$ be a subgroup of 
$GL(2,\mathbb C).$ 
$G$ is said to be extendable if 

$$ 0\rightarrow \{\pm Id\} \rightarrow G\stackrel{\pi}{\rightarrow} 
K\rightarrow 0$$
\noindent is exact, and $G$ is said to be non-extendable 
if $\pi_{\mid G}\, :G\rightarrow K$ 
is an isomorphism.
\end{definition}

Let $B=\{ (x_{0},x_{1})\in \mathbb P^{1}\, \mid \, \beta (x_{0},x_{1})=0\}.$
\begin{remark} \label{marina}
Let $G\subset{\rm{Aut}}(C)$ and $K=\pi(G).$ It holds: 
$1)$ $G$ is extendable iff $i\in G;$ 
$2)$ $B$ is a $K$-invariant. 
\end{remark}
Obvious.\qed
\vspace{.10in}\noindent

The finite subgroups $K\subset\mathbb PGL(2, \mathbb C)$ are well known 
\cite{Kl}, and it is easy to find 
the $K$-invariant polynomials 
$\beta$. 

Then we will show that 
the data $K$ and $\beta$ 
{\underline{uniquely determine an extendable group}} 
$G$ acting on $C$.

\begin{proposition}\label{origene}
Let $K\subset\mathbb PGL(2,\mathbb C)$ be a finite subgroup and 
$B=\{ \beta (x_{0},x_{1})=0\}$ a $K$-invariant reduced divisor of degree $6$. 
There exists a {\rm{unique}} $G\subset GL(2,\mathbb C)$ such that 
$i)$ $G$ is a group of automorphisms of 
$\frac {\mathbb C[x_{0},\, x_{1},\, z]}
{(z^{2}-\beta (x_{0},\, x_{1}))}$ 
and $ii)$ $G$ is extendable.
\end{proposition}
\proof {\it{Unicity}}. If 
$G_{1}, G_{2}\subset GL(2)$ 
satisfy 
the claim and 
$h_{1}\in G_{1}-G_{2},$ there exists 
$h_{2}\in G_{2}$ such that $\pi h_{1}=\pi h_{2}$, 
since $K=\pi G_{1}=\pi G_{2}$. But 
$i\in G_{1}\cap G_{2},$ then $h_{1}=ih_{2}\in G_{2}$: 
a contradiction.

\noindent
{\it{Existence.}} 
By (\ref{pgln}) we have:  

\begin{equation}\label{induced}
0\rightarrow \{\pm Id\} \rightarrow 
\widehat{K}\stackrel{\pi_{1}}{\rightarrow}K \rightarrow 0.
\end{equation}
\noindent
Using the $K$-invariance of $B$, {\it{define}} $\lambda \, :\widehat{K}\rightarrow \bbbc^{\star}$ by
$(\pi(\widehat {k}))^{\star}(\beta )=\lambda (\widehat {k})\beta.$ Since 
${\rm{deg}}(\beta )$ is even, $\lambda$ descends to the quotient 
$K=\frac{\widehat{K}}{\{\pm Id\}}$ and we will not distinguish this latter character from $\lambda$. {\it{Define:}}
\vspace{.10in}\par\noindent
\begin{equation}\label{definizione}
G=\{\pm \sqrt{(\lambda (k))^{-1}}\,\widehat{k}\, \mid \, \widehat {k}\in \widehat {K} \, ,\, k=\pi(\widehat {k})\}.
\end{equation}
\vspace{.10in}\noindent
$G$ is the claimed group. Notice that $G\subset GL(2,\mathbb C)$ and $G$ 
is a group. We can write the $G$-action over ${\cal{R}}$:

\begin{equation}\left\{ \begin{array}{cc}
g(x_{0})=& \pm \sqrt{(\lambda (k))^{-1}} k(x_{0})\\
g(x_{1})=& \pm \sqrt{(\lambda (k))^{-1}} k(x_{1})\\
g(z)=& (\lambda (k))^{-1}z= det(g)z\\ 
\end{array} \right. 
\end{equation}
\vspace{.10in}\par\noindent
where $g=\pm \sqrt{(\lambda (k))^{-1}}\,\widehat{k}.$
It remains to prove that  
$$0\rightarrow \{\pm Id\} 
\rightarrow G\stackrel{\pi}{\rightarrow}K \rightarrow 0$$\noindent is exact. 
By definition $\pi$ is surjective. Let $g\in \, {\rm{ker}} \pi$. 
By (\ref{definizione}) there exists $\widehat {k}\in\widehat{K}$ such that 
$g=\pm \sqrt{(\lambda (k))^{-1}}\,\widehat{k}$. Then $\widehat {k}=\pm Id$; 
in particular $\lambda(k)=1$ and it follows $g=\pm Id.$\qed\vspace{.10in}
\noindent

Actually we have shown:
\begin{theorem}\label{realta} 
The couples $(C,\, G)$ where $C$ is a fixed curve of genus $2$ and 
$G\subset {\rm{Aut}}(C)$ is extendable are in bijection with the classes 
$(K,\, B)$ up to ${\rm{Aut}(}\mathbb P^{1})$ where 
$K\subset\mathbb PGL(2,\mathbb C)$ is a finite subgroup and 
$B$ is a degree $6$, $K$-invariant, reduced divisor.
\end{theorem}\vspace{.10in}\noindent

We want to classify non-extendable groups. 
If $G$ is non-extendable we need to understand how $G$ fits into 
$G^{'}={\rm{Aut}}(C)$. To this end, set $K^{'}=\pi(G^{'})$, 
and restrict (\ref{gln}) to 
$G^{'}$:
\begin{equation}\label{restri}
 0\rightarrow \{\pm Id\} \rightarrow G^{'}\stackrel{\pi}{\rightarrow} K^{'}\rightarrow 0.
\end{equation}
\noindent Notice that in general 
it is {\it{not}} true that $G\subset SL(2,\mathbb C).$ We like 
to consider subgroups $K\subset K^{'}$ and their 
liftings to $G^{'}$.
\begin{definition}\label{splitta}
A subgroup $K\subset K^{'}$ is said to be of  splitting-type if 
$(\pi)^{-1}(K)$ is splitted 
(i.e. $(\pi)^{-1}(K)=K\times \{ \pm Id \}.$ 
Otherwise it is of non splitting type.
\end{definition}
Obviously we have the following remarks that we write for further reference:
\begin{remark}\label{melitone}
$K$ is splitting if and only if there exist a lifting 
$\epsilon :K \rightarrow G^{'}$ such that 
$\pi \circ \epsilon =Id_{K}.$ In particular $K$ is splitting if and only if there exists a nontrivial 
homomorphism $\epsilon: K\stackrel {}{\rightarrow} \{\pm Id\}.$
\end{remark}

\begin{remark}
If $G$ in $(C\, , \, G)$ is non-extendable then 
$K=\pi (G)$ is of splitting type.
\end{remark}
On the other hand if If $G$ is extendable both cases for $K=\pi(G)$ 
may occur, 
but the following 
case is easy to describe:
\begin{remark}\label{ippolito}
If $G$ in $(C\, ,\, G)$ is extendable and $\pi(G )=K$ is of 
splitting type then $G=K\times \{\pm Id\}.$
\end{remark}
The following corollary gives the analogue of \ref{origene} for the 
non-extendable groups:
\begin{corollary}\label{nonsuper}
Let $K$ and $B$ as in \ref{origene}. The set of the 
couples $(C,\, ,G)$ where $G$ is non-extendable 
is in bijection with the set of the liftings 
$\epsilon:K\rightarrow G_{s}$ where 
$G_{s}$ is the unique extendable group constructed in 
the proof of \ref{origene} through the data $K$ and $B$.  
\end{corollary}
\proof Trivial. \qed\vspace{.10in}

If $G$ 
is of splitting type there exists a lifting $K\stackrel{\epsilon}{\rightarrow}
GL(2,\mathbb C)$ and let us denote by 
$\mu:K\stackrel{}{\rightarrow}\mathbb C^{\star}$ 
the character uniquely defined by the relation 
$\mu (k)\beta=(\epsilon (k))^{\star}( \beta ).$

\begin{proposition}\label{gregorio}
The isomorphism class $(K\, ,\, B)$ 
up to ${\rm{Aut}(}\mathbb P^{1})$ 
induces $(C\, , \, G_{s})$ where $G_s$ is splitted if and only if 
there exists $i)$ a lifting 
$K\stackrel{\epsilon}{\rightarrow}GL(2, \mathbb C)$ and $ii)$ 
a square root $\nu$ of the character $\mu$ associated to 
$\epsilon$ and $\beta.$
\end{proposition}
\proof
Assume that the procedure described in \ref{origene} gives a splitted couple 
$(C\, ,\, G_{s}).$
By \ref{ippolito}, $i)$ follows. From the proof of \ref{origene} and by 
$\epsilon: K\stackrel{}{\rightarrow}GL(2, \mathbb C)$ we have that in
(\ref{induced}) 
\begin{equation}\label{binduced}
\widehat {K}=
\{\widehat {k}= 
\pm \sqrt{(det 
(\epsilon (k)))^{-1}}\epsilon (k)\, \mid \, k\in K\, \}.
\end{equation}
\noindent
Moreover the character $\lambda$ to construct the claimed 
$G_{s}$ is by definition 
$$\lambda (k)=(det (\epsilon (k)))^{-3}\mu (k).$$
\noindent By \ref{origene} 
$$G_{s}=\{ \pm \sqrt{(det (\lambda (k)))^{-1}}\,{\widehat {k}}\, \mid \, {\widehat {k}}\in{\widehat {K}} \},$$ \noindent
then by the form  ${\widehat {k}}$ in (\ref{binduced}) it holds 
$$\pm \sqrt
{ (\lambda (k))^{-1}}
\sqrt{(det (\epsilon (k)))^{-1}}{\epsilon (k)}=
\pm \sqrt{((det (\epsilon (k)))^{-3}\mu(k))^{-1}}\sqrt{(det (\epsilon (k)))^{-1}}{\epsilon (k)}$$\noindent that is, if we set $\rho(k)=\sqrt{\lambda(k)(det (\epsilon (k))^{-1}}$ we can write the claimed square root 
$\nu(k)=\frac{det(\epsilon(k))}{\rho(k)}.$

\noindent
{\it{Viceversa.}} Assume that $i)$ and $ii)$ hold. We {\it{define}} 
$\rho(k)=\frac{det(\epsilon(k))}{\nu(k)}$ and the same computation 
in reverse order shows that $(\epsilon)^{'}=\rho \epsilon$ is a lifting  
${\epsilon}^{'}\, :K\rightarrow G_{s}$. The by \ref{melitone} 
$G_{s}$ is splitted.\qed\vspace{.10in}


{\begin{center}{\bf{Bolza classification}}\end{center}}
In \cite{Kl} Klein shows the finite groups acting on $\mathbb P^{1}$:
{\samepage $$\begin{tabular}{|l|c|l|}\hline
Group K & order special orbits  & Order of K\\ \hline
$\mathbb Z/n$ & 1\, , \, 1 & n \, \, Cyclic\\ \hline
${\cal {D}}_{n}$ & n\, ,\, n\, ,\, 2 & 2n\, \, Dihedral\\ \hline
${\cal {A}}_{4}$ & 6\, ,\,4\, ,\, 4 & 12\, \, tetraedral\\ \hline
${\cal {S}}_{4}$ & 12\, ,\, 8\, ,\,6 & 24\, \, esaedral o octaedral\\ \hline
${\cal {A}}_{5}$ & 30\, ,\,30\, ,\,12 & 60 \,\, 
Icosaedral o dodecaedral\\ \hline
\end{tabular}$$
\vbox to 0pt{$${\rm{(Table \,\,I}})$$}}
\vspace{.10in}\par\noindent
In our case $B$ is a reduced $K$-invariant divisor then ${\cal {A}}_{5}$ does not occur In the same book we find the groups 
$\widehat {K}$ such that  $0\rightarrow \{\pm Id\}\rightarrow\widehat {K}\rightarrow K\rightarrow 0$ is exact:

$$\begin{tabular}{|l|c|l|}\hline
Group $\widehat{K}$ & Generator  & Relations\\ \hline
$\mathbb Z/2n$ & $\zeta =\left( \begin{array}{cc}
e^{i\pi /n}& 0\\
0& e^{-i\pi/n}
\end{array}\right) $&  \, \, $\langle \zeta\mid {\zeta}^{2n}=1\rangle$\\ \hline
${\cal {D}}_{n}\mathbb Z/2$ &$\zeta =\left( \begin{array}{cc}
e^{i\pi /n}& 0\\
0& e^{-i\pi/n}
\end{array}\right) $, $\eta =\left( \begin{array}{cc}
0& i\\
i& 0
\end{array}\right) $ & $\langle \zeta\, ,\, \eta \mid \begin{array}{c}
{\zeta}^{2n}={\eta}^{4}=1\\
{\zeta}^n={\eta}^2\\{\eta}^2\zeta=\zeta {\eta}^2
\end{array} 
\rangle$
\\ \hline
${\widehat{\cal {A}}}_{4}$ & $\begin{array}{c}
\zeta =1/2(i-1)\left( \begin{array}{cc}
1& -i\\
1& i
\end{array}\right)\\ \eta =(i+1)\left( \begin{array}{cc}
i& -i\\
-1& -1
\end{array}\right) \end{array} $ & $\langle\ \zeta\, ,\, \eta \mid \begin{array}{l}
{\zeta}^{3}={\eta}^{3}={{\eta}{\zeta}}^{4}=1\\
({\zeta}{\eta})^2=({\eta}{\zeta})^2\\ ({\eta}{\zeta})^{2}{\zeta}={\zeta} ({\eta}{\zeta})^2\\({\eta}{\zeta})^{2}{\eta}={\eta} ({\eta}{\zeta})^2

\end{array}
\rangle$ \\ \hline
${\widehat{\cal {S}}}_{4}$ & $\begin{array}{c}
\zeta =1/2(i-1)\left( \begin{array}{cc}
1& -i\\
1& i
\end{array}\right)\\ \eta =\frac{1}{\sqrt{2}}\left( \begin{array}{cc}
i+1& 0\\
0& 1-i
\end{array}\right) \end{array} $ &  $\langle\ \zeta\, ,\, \eta \mid \begin{array}{l}
{\zeta}^{3}={\eta}^{8}={{\eta}{\zeta}}^{4}=1\\
{\eta}^4{\zeta}={\zeta}{\eta})^4\\

\end{array}
\rangle$ \\ \hline
\end{tabular}$$
$${\rm{(Table II)}}$$
\vspace{.10in}\par\noindent
If one wants to look directly to the quoted book, notice that
if  ${\widehat {K}}={\widehat {\cal {A}}}_{4}$ 
then, in the book's notation, we have:

$$\zeta\eta=\left( \begin{array}{cc}
i&0\\0&-i
\end{array} \right)\, \, , \,\,
\eta \zeta=\left( \begin{array}{cc}
0&1\\-1&0
\end{array} \right)$$
\vspace{.10in}\par\noindent
while if ${\widehat {K}}={\widehat {\cal {S}}}_{4}$ 
then 
$${\eta}^{2}=\left( \begin{array}{cc}
i&0\\0&-i
\end{array} \right)\, \, , \,\,
(\eta \zeta)^{2}=\left( \begin{array}{cc}
0&-i\\i&0
\end{array} \right)
\, \, , \,\,
(\zeta \eta)^{2}=\left( \begin{array}{cc}
0&-1\\1&0
\end{array} \right).$$
\vspace{.10in}\par\noindent
 We have just noticed that the case ${\cal {A}}_{5}$ does not occur but it is not the unique one:

\begin{lemma}\label{tabellatret}

The groups of table $(I)$ which give a 
couple $(C\, ,\, G)$ where $G$ is extendable
 are completely 
classified in the following table:

 $$\begin{tabular}{|l|c|l|}\hline
 K &  G  & C\\ \hline
${\bbbz}/6$ & ${\bbbz}/6\times {\bbbz}/2$ & $z^{2}=x_{1}^{6}-x_{0}^{6}$\\ \hline

${\bbbz}/5$ & ${\bbbz}/10$ & $z^{2}=x_{0}(x_{1}^{5}-x_{0}^{5})$\\ \hline

${\bbbz}/4$ & ${\bbbz}/8$ & $z^{2}=x_{1}x_{0}(x_{1}^{4}-x_{0}^{4})$\\ \hline

${\bbbz}/3$ & ${\bbbz}/6$ & $z^{2}=(x_{1}^{3}-x_{0}^{3})(x_{1}^{3}+x_{0}^{3})$\\ \hline

${\bbbz}/2$ & ${\bbbz}/4$ & $z^{2}=(x_{1}^{2}-x_{0}^{2})(x_{1}^{2}-4x_{0}^{2})(x_{1}^{2}-9x_{0}^{2})$\\ \hline

${\bbbz}/2$ & ${\bbbz}/2\times {\bbbz}/2$ & $z^{2}=x_{1}x_{0}(x_{1}^{2}-x_{0}^{2})(x_{1}^{2}-4x_{0}^{2})$\\ \hline

$\langle{\rm{id}}\rangle$ & ${\bbbz}/2$ & 
$z^{2}={\rm{generic\,\, polynomial\,\, of \,\,degree\,\, 6}}$\\ \hline

${\cal {D}}_{6}$ & $\langle T\, U\,\mid\,\begin{array}{c}
T^{2}=U^{6}=(TU)^{4}=1\\ T(TU)^2=TU^{2}T\\ U(TU)^2=(TU)^{2}U \end{array}\rangle  $& $z^{2}=(x_{1}^{6}-x_{0}^{6})$\\ \hline 

${\cal {D}}_{4}$ & $\langle T\, U\,\mid\,\begin{array}{c}
T^{2}=U^{8}=(UT)^{4}=1\\ TU^4=U^{4}T \end{array} \rangle  $& $z^{2}=x_{0}x_{1}(x_{1}^{4}-x_{0}^{4})$\\ \hline

${\cal {D}}_{3}$ & ${\cal {D}}_{6}$ & $z^{2}=(x_{0}^{3}-2x_{1}^{3})(x_{1}^{3}-2x_{0}^{3})$\\ \hline

${\cal {D}}_{3}$ & $\langle T\, U\,\mid\,\begin{array}{c}
T^{3}=U^{4}=(TU)^{4}=1\\ (TU)^2=U^{2}\\ TU^2=U^{2}T \end{array} \rangle $& $z^{2}=(x_{1}^{6}-x_{0}^{6})$\\ \hline 

${\cal {D}}_{2}$ & ${\bbbz}/4\times {\bbbz}/2$& $z^{2}=x_{0}x_{1}(x_{1}^{2}-4x_{0}^{2})(x_{0}^{2}-4x_{1}^{2})$\\ \hline

${\cal {D}}_{2}$ & ${{\cal{H}}}$& $z^{2}=x_{0}x_{1}(x_{0}^{4}-x_{1}^{4})$\\ \hline

${\cal {A}}_{4}$ & ${\widehat{\cal {A}}}_{4}$ & $z^{2}=x_{0}x_{1}(x_{0}^{4}-x_{1}^{4})$\\ \hline 

$ {\cal {S}}_{4}$ & $\langle T\, U\,\mid\,\begin{array}{c}
T^{3}=U^{8}=(UT)^{2}=1\\ TU^4=U^{4}T \end{array} \rangle  $ & $z^{2}=x_{0}x_{1}(x_{0}^{4}-x_{1}^{4})$\\ \hline 
\end{tabular}$$
\vbox to 0pt{$${\rm{Table \,(III)}}$$}
\end{lemma}

\noindent
\proof {\it{First step: to find the occurrences of $\beta$}}. 
This is easy since for each $K$ in Table $(I)$ we have to find 
which unions of orbits have order $6$.\vspace{.10in}

\noindent
{\it{Second step: to find $G$.}} By the first step 
we know $B$ and $K$. The procedure described in 
\ref{origene} gives the result. For example, 
we show how to obtain the extendable group associated to 
${\cal {S}}_{4}.$ In particular we will see that it is 
{\it{different}} from $\widehat{{\cal {S}}}_{4}.$ 
Let us consider the Table $(II)$. 
We find two generators of $\widehat{{\cal {S}}}_{4}:$
$$\begin{array}{cc}
\zeta =1/2(i-1)\left( \begin{array}{cc}
1& -i\\
1& i
\end{array}\right)\\ \eta =\frac{1}{\sqrt{2}}\left( \begin{array}{cc}
i+1& 0\\
0& 1-i
\end{array}\right). \end{array}$$
\vspace{.10in}
\noindent
The character $\lambda$, 
to start with \ref{origene}, 
can be obtained by the following game:

$$\begin{array}{c}
\lambda (\zeta)\beta=(\zeta)^{\star}(\beta)=\beta\\
\lambda (\eta)\beta=(\eta)^{\star}(\beta)=-\beta.
\end{array}$$
\vspace{.10in}\par\noindent
By construction, the solution is 
$G=\left\{ \pm\sqrt{(\lambda(k))^{-1}}\,{\widehat{k}}\, 
\mid\, {\widehat{k}}\in {\widehat{{\cal {S}}}}_{4}\,\right \}$ and in 
$G$ there are:

$$\begin{array}{cc}
T=(\zeta)^{2} =i/2(i-1)\left( \begin{array}{cc}
1& 1\\
i& -i
\end{array}\right)\\ U=i\eta =\frac{1}{\sqrt{2}}\left( \begin{array}{cc}
i-1& 0\\
0& i+1
\end{array}\right). \end{array}$$
\vspace{.10in}\par\noindent
Let $H=\langle\, T\, ,\,U\,\rangle$ be the subgroup generated by $T$ and $U$. We want to show that $H=G$. Obviously: $H\subset G$ and $U^{4}=-Id.$ 
Then we can restrict  (\ref{gln}) to $H$ and it gives 
$K_{1}\subset \mathbb PGL(2,\mathbb C)$ such that 
$$0\rightarrow \langle\, U^{4}\,\rangle\rightarrow H\stackrel{\pi}{\rightarrow}K_{1}\rightarrow0$$\noindent
is exact. The task is to show that $K_{1}={\cal {S}}_{4}$. 
That is, we have to show that ${\cal {S}}_{4}\subset K_{1}.$ 
Set $u=\pi(U)$ and $t=\pi(T).$ 
From Table $(II)$ we have $u=\pi(\eta)$ and 
$t=\pi({\zeta}^{2})$. We conclude by (\ref{pgln}) restricted to 
${\cal {S}}_{4}.$\qed\vspace{.10in}\par\noindent

\begin{remark}\label{tablequattro}
The way we prove \ref{tabellatret} gives an explicit description of 
the extendable groups as subgroups of $GL (2, \mathbb C)$. 
They are listed in Table $(IV)$, see appendix.  
\end{remark}

To end the classification we have to find which groups in Table $(III)$ 
are splitted and in the affirmative case to classify all the liftings 
$\epsilon: K\rightarrow G$ {\underline{such that}} 
$C/\epsilon(K)=\mathbb P^{1}$. It requires only a few 
basic facts on curves theory: essentially that 
if $C$ is a genus $2$ curve, $G\subset{\rm{Aut}}(C)$ and 
$C/G$ is elliptic then $G={\bbbz}/2;$ this follows from 
the Hurwitz formula and the easy monodromy argument that an Abelian 
covering over an elliptic curve has at least two branch points. 
In some cases to find the splitted group $G_{s}$, instead of 
\ref{gregorio}, we can use a more direct argument.

\begin{lemma}\label{simm}
If $K={\cal {A}}_{4}$ or $K={\cal {S}}_{4}$ then $K$ is non splitting.
\end{lemma}
\proof
The proofs are similar. We only do the case $K={\cal {S}}_{4}.$ If  
${\cal {S}}_{4}$ were splitting, then the corresponding extendable  
group in Table $(III)$ would be $G={\cal {S}}_{4}\times {\bbbz}_{2};$ 
a contradiction, because in $G$ there is an element of order 
$8$.\qed
\vspace{.10in}

\noindent
{\begin{center}{\bf{The dihedral case}}\end{center}}
If $K={\cal {D}}_{n}$ we like to consider 
two cases depending on 
the parity of $n$.
\begin{lemma}\label{celso}
Let $K={\cal {D}}_{n}$, then $K$ is splitting if and only if $n$ is odd.
\end{lemma}
\proof ${\cal {D}}_{n}\subset \mathbb PSL(2)$ is given by:
$$\langle\, \left[ \left( \begin{array}{cc} 0&1\\1&0 \end{array} \right) \right] \,,\, \left[ \left( \begin{array}{cc} 1&0\\0&\xi \end{array} \right) \right] \,\rangle.$$
\noindent
If ${\cal {D}}_{n}\stackrel{\epsilon}{\rightarrow}GL(2)$ is a lifting, 
the preimage of 
$\left[ \left( \begin{array}{cc} 0&1\\1&0 \end{array} \right)\right]$ 
is $\{\pm \left( \begin{array}{cc}0&1\\1&0 \end{array} \right)\}$ 
and every liftings of $\left[ \left( \begin{array}{cc} 1&0\\0&\xi \end{array} \right) \right]$ has the following form: 
${\xi}^{i}\left( \begin{array}{cc} 1&0\\0&\xi \end{array} \right).$ 
By definition, 
$\epsilon$ is an isomorphism over its image then the relation defining 
$\epsilon ({\cal {D}}_{n})$ forces
$${\xi}^{i} \left( \begin{array}{cc} \xi &0\\0&1 \end{array} \right) =
({\xi}^{i} \left( \begin{array}{cc} 1&0\\0&\xi \end{array} \right))^{-1}$$
\noindent
to hold. Then ${\xi}^{2i+1}=1$, which 
has a solution if and only if $n$ is odd.\qed
\vspace{.10in}

By \ref{celso} and by Table $(III)$ the case 
$K={\cal {D}}_{n}$ splitting is achieved 
applying \ref{gregorio} to $K={\cal {D}}_{3}.$ We recall: 
\begin{remark}\label{ilario}
If $n$ is odd then ${\bbbz}_{2}$ is the group of 
${\cal {D}}_{n}$-linear characters.
\end{remark}

Let ${\cal {D}}_{3}\stackrel{\epsilon}{\rightarrow}GL(2)$ be a lifting, 
for example:
$$\epsilon \left[ \left( \begin{array}{cc} 1&0\\0&\xi \end{array} \right) \right]= \left( \begin{array}{cc} 1&0\\0&\xi \end{array} \right) \, ,\, 
\epsilon \left[ \left( \begin{array}{cc} 0&1\\1&0\ \end{array} \right) \right]= \left( \begin{array}{cc} 0&1\\1&0 \end{array} \right)$$
\noindent
where $\xi=e^{2i\pi /3}.$ There are two cases. If $\beta=x_{1}^{6}-x_{0}^{6}$ the character $\mu$ induced by $\epsilon$, $\beta$ is:
$$\mu \left( \begin{array}{cc} 1&0\\ 0&\xi \end{array} \right)=1 \, ,\, 
\mu \left( \begin{array}{cc} 0&1\\1&0 \end{array} \right)=-1.$$
\noindent
By \ref{ilario}, 
it does not exist a character 
$\nu$ such that ${\nu}^{2}={\mu}.$
\vspace{.10in}\noindent
In the other case 
$\beta= (x_{0}^{3}-2x_{1}^{3})(x_{1}^{3}-2x_{0}^{3}),$ and we have 
$$\mu \left( \begin{array}{cc} 1&0\\0&\xi \end{array} \right)=1 \, ,\,
\mu \left( \begin{array}{cc} 0&1\\1&0 \end{array} \right)=1.$$
\noindent
We easily see that $\nu=1$ satisfy the condition of \ref{gregorio}. 
If ${\epsilon}_{1}=det(\epsilon)\epsilon$, 
${\epsilon}_{1}:{\cal {D}}_{3}\stackrel{}{\rightarrow}{\cal {D}}_{6}\subset GL(2)$ is a lifting and 
$${\epsilon}_{1}({\cal{D}}_{3})=
G_{1}=\langle\, \left( \begin{array}{cc} \xi &0\\0&{\xi}^{2} \end{array} \right)\, ,\,\left( \begin{array}{cc} 0&-1\\-1&0 \end{array} \right) \, 
\rangle.$$

To conclude the case ${\cal{D}}_{3}$ we need:
\begin{remark}\label{pistol}
Let $K\subset \mathbb PGL(2, \mathbb C)$ with a fixed lifting 
$K\stackrel{{\epsilon}_{1}}{\rightarrow}GL(2, \mathbb C)$. 
Then for every lifting 
${\epsilon}_{2}:K\stackrel{}{\rightarrow}GL(2, \mathbb C)$ it holds 
${\epsilon}_{2}={\rho}{\epsilon}_{1}$ where $\rho$ is a character of $K$.
\end{remark}
\vspace{.10in}
Then by \ref{pistol} and \ref{ilario}, we have another lifting 
${\epsilon}_{2}={\rho}{\epsilon}_{1},$ with 
    
$$\rho \left[ \left( \begin{array}{cc} 1&0\\0&\xi \end{array} \right)\right]=1 \, ,\,
\rho \left[ \left( \begin{array}{cc} 0&1\\1&0 \end{array} \right)\right]=-1,$$
\noindent
and
$${\epsilon}_{2}({\cal {D}}_{3})=
G_{2}=\langle\, \left( \begin{array}{cc} \xi &0\\0&{\xi}^{2} \end{array} \right)\, ,\,\left( \begin{array}{cc} 0&1\\1&0 \end{array} \right) \, \rangle.$$
\noindent
We sum up the dihedral case in the following lemma: 

\begin{lemma}\label{pucia}
If $K={\cal {D}}_{n}$ is splitting then $K={\cal {D}}_{3},$ 
$\beta=(x_{0}^{3}-2x_{1}^{3})(x_{1}^{3}-2x_{0}^{3})$ 
and there are two liftings 
${\epsilon}_{i}:{\cal {D}}_{3}
\rightarrow{\cal {D}}_{6}\subset 
GL(2, \mathbb C)$ where $i=1,2$ such that if 
$G_{i}={\epsilon}_{i}({\cal {D}}_{3})$ then $C/G_{i}=I\!\!P^{1}.$
\end{lemma}
\vspace{.10in}\noindent

\begin{center}{\bf{The cyclic case}}\end{center} 
Even in the case $K=\mathbb Z_{n}$ the behaviour depends on the parity.

\begin{lemma}\label{amore}
Let $K=\mathbb Z_{n}$ with $n$ odd. If $G$ is the corresponding 
extendable group then 
$G=\mathbb Z_{2n},$ and it is splitted. Moreover every lifting 
$\epsilon: K\rightarrow G$ gives {\rm{the same}} 
subgroup of $G$ and $C/\epsilon(K)=\mathbb P^{1}.$
\end{lemma}
\proof An easy computation with the cases 
$n=5$, $n=3$ in Table $(III)$.\qed\vspace{.10in}

If $n$ is even we have more cases.
\begin{lemma} If $n$ is even, $K=\mathbb Z_{n}$ 
and $B$ contains some special orbits then the corresponding 
extendable group is non-splitted.
\end{lemma}
\proof
The special orbit in $B=\left\{ \beta =0\right\}$ is given by 
$x_{0}x_{1}=0$.
By table $(I)$ and table $(III)$ we have to consider only the case with 
$n=4$ or  $n=2.$ The cyclic subgroup of 
$I\!\!PGL(2, \mathbb C)$ is $\langle\,     
\left[ \left( \begin{array}{cc} 1&0\\0&\xi \end{array} \right)\right]
\rangle,$ where $\xi=e^{2i\pi/n}.$ Let 
$\epsilon \left[ \left( \begin{array}{cc} 1&0\\0&\xi \end{array} \right)\right]= \left( \begin{array}{cc} 1&0\\0&\xi \end{array} \right)$ be a lifting. 
Using the notation of \ref{gregorio}, $\mu \left[ \left( \begin{array}{cc} 1&0\\0&\xi \end{array} \right)\right]=\xi$. Then it 
does not exists a character $\nu$ di ${\bbbz}_{n},$ 
such that ${\nu}^{2}=\mu$. By 
\ref{gregorio} we conclude.\qed\vspace{.10in}\noindent

By table $(III)$ we have only two cases to consider: 
$(K\, ,\, B)= (\mathbb Z_{6}\, ,\, x_{1}^{6}-x_{0}^{6}))$ and  
$(K\, ,\, B)= 
(\mathbb Z_{2}\, ,\, x_{1}^{2}-x_{0}^{2})(x_{1}^{2}-
x_{0}^{2})(x_{1}^{2}-x_{0}^{2})).$

\begin{lemma}\label{facile}
Let $(K\, ,\, B)= (\mathbb Z_{6}\, ,\, x_{1}^{6}-x_{0}^{6}))$ 
then $K$ is splitting and it has two liftings 
$\epsilon_{i}:
{\bbbz}_{6}\stackrel{}{\rightarrow}{\bbbz}_{12}\subset 
GL(2, \mathbb C)$, $i=1,2.$ 
Letting $G_{i}={\epsilon}_{i}({\bbbz}_{6})$ it holds that 
$C/G_{i}=\mathbb P^{1}$ and 
$$G_{1}=\langle \left( \begin{array}{cc} e^{2i\pi/3}&0\\0&e^{i\pi/3} \end{array} \right)\rangle,\,\,\, G_{2}=\langle \left( \begin{array}{cc} -e^{2\pi/3}&0\\0& -e^{i\pi/3} \end{array} \right)\rangle.$$
\end{lemma}
\proof An easy computation.\qed \vspace{.10in}\par\noindent
\begin{lemma}\label{simplicio}
Let $(K\, ,\, B)= 
(\mathbb Z_{2}\, ,\, x_{1}^{2}-x_{0}^{2})(x_{1}^{2}-x_{0}^{2})(x_{1}^{2}-x_{0}^{2}))$. Then $K$ is splitting and it has two liftings 
${\epsilon}_{i}:\mathbb Z_{2}\stackrel{}{\rightarrow}{\bbbz}_{2}\times {\bbbz}_{2}\subset GL(2)$, $i=1,2.$ Moreover if  
$G_{i}={\epsilon}_{i}({\bbbz}_{2})$ then $C/G_{i}$ is an elliptic curve.
\end{lemma}
\proof
The group 
$\mathbb Z_{2}\times \mathbb Z_{2}\subset GL(2\mathbb Z,)$ 
is 
$G=\langle \left( \begin{array}{cc} 1&0\\0&-1 \end{array} \right)\, ,\, \left( \begin{array}{cc} -1&0\\0&1 \end{array} \right)\,\rangle.$ 
$K$ is generated by $\left[\left( \begin{array}{cc} 1&0\\0&-1 \end{array} \right) \right],$ and its two liftings are 
$\pm\left( \begin{array}{cc} 1&0\\0&-1 \end{array} \right).$ Finally if 
$H_{0}=\langle \left( \begin{array}{cc} 1&0\\0&-1 \end{array} \right)\rangle$ 
and 
$H_{1}=\langle\left( \begin{array}{cc} -1&0\\0&1 \end{array} \right)\,\rangle$ 
then $x_{0},$ $x_{1}$ in canonical ring correspond 
to the $1$-form invariant 
by $H_{0}$ and respectively, by $H_{1}$. In particular 
$C/H_{i}$ has genus $1$ where $i=0,1.$\qed\vspace{.10in}\noindent

The classification result is
\begin{theorem}\label{trefolds}
Let $C$ be a curve of 
genus $2$ and let 
$G\subset{\rm{Aut}}(C)$ be a non trivial subgroup such that
$C/G=\mathbb P^{1}.$ There are only $21$ types of couples
 $(C\, ,\, G)$. Moreover 
$15$ have extendable type and they are listed in table $(III)$. 
The remaining $6$ types are listed below:
$$\begin{tabular}{|l|c|l|}\hline
 $K\sim\epsilon(K)=G$ & $G_{s}$ &generators of $\epsilon(K)$\\ \hline

$\mathbb Z/3$ &$\mathbb Z/3 \times \mathbb Z/2$&
$\left( \begin{array}{cc}
e^{2i\pi/3}
& 0\\
0& e^{4i\pi/3}
\end{array}\right) $\\ \hline

$\mathbb Z/5$ &$\mathbb Z/5 \times \mathbb Z/2$&
$\left( \begin{array}{cc}
e^{2i\pi/5}
& 0\\
0& e^{4i\pi/5}
\end{array}\right) $\\ \hline

$\mathbb Z/6$ &$\mathbb Z/6 \times \mathbb Z/2$&
$\left( \begin{array}{cc}
e^{2i\pi/3}
& 0\\
0& e^{i\pi/3}
\end{array}\right) $\\ \hline
$\mathbb Z/6$ &$\mathbb Z/6 \times\mathbb Z /2$ &
$\left( \begin{array}{cc}
-e^{2i\pi /3}& 0\\
0& -e^{i\pi/3}
\end{array}\right) $\\ \hline
${\cal {D}}_{3}$ & ${\cal {D}}_{6}$ & $\left( \begin{array}{cc}
0&-1\\
-1& 0
\end{array}\right)$,
$\left( \begin{array}{cc}
e^{2i\pi /3}& 0\\
0& e^{i\pi/3}
\end{array}\right) $\\ \hline
${\cal {D}}_{3}$ & ${\cal {D}}_{6}$ & $\left( \begin{array}{cc}
0&1\\
1& 0
\end{array}\right)$,
$\left( \begin{array}{cc}
-e^{2i\pi /3}& 0\\
0&- e^{i\pi/3}
\end{array}\right) $\\ \hline 
\end{tabular}$$
\end{theorem}
\proof It follows from \ref{tabellatret}, 
\ref{simm}, 
\ref{pucia}, 
\ref{amore},
\ref{facile},
\ref{simplicio}.\qed\vspace{.10in}\noindent

The full classification of the case with $p_{g}=q=2$ and $S$ 
not of general type would require to compute all the unramified $G$ actions 
$C_{1}\rightarrow C_{1}/G$ where $g(C_{1}/G)=2$, for each occurrence of $G$ mentioned in \ref{trefolds}. We think that the outcome is not worthy of the 
effort. 
However since $S$ is GH,  by \cite[Theorem B]{Ca1}, 
\cite[Theorem C]{Ca1} then we can say:
\begin{theorem}\label{finiamola}
Each irreducible component of the moduli space of surfaces 
with $p_{g}=q=2$ and not of Albanese general type is given by 
$\cM(\Pi,4)$, the moduli space of surfaces isogenus to a product 
with fundamental group 
$\Pi$, Euler number $4$ and each component is specified by a 
fixed isomorphism $\Pi(S)\rightarrow \Pi$ which fits into the exact sequence
$$
0\rightarrow \Pi_{1}(C_{1})
\times \Pi_{1}(C_{2})\rightarrow \Pi\rightarrow G\rightarrow 0
$$
\noindent
such that the factors $\Pi_{1}(C_{i})$ are normal in $\Pi$, 
the orbifold exact sequences of the coverings $C_{1}\rightarrow C=C_{1}/G$, 
$C_{2}\rightarrow \mathbb P^{1}=C_{2}/G$, 
$$
0\rightarrow \Pi_{1}(C_{i})\rightarrow \Pi(i)\rightarrow G\rightarrow 0
$$\noindent are such that there is 
no element of $\Pi$ mapping in each $\Pi(i)$ to an element of finite order 
and $G$ embeds in $\rm{Out}(\Pi(C_{i}))$ by the above 
sequence where $i=1,2$ and $G$ is one of the groups classified in 
\ref{trefolds}.
\end{theorem}

\noindent
Zucconi Francesco\\
Universit\`{a} di Udine Dipartimento di Matematica e Informatica\\
Via delle Scienze 206 33100 Udine, Italia.\\
e-mail zucconi@dimi.uniud.it

\end{document}